\documentclass[dvips,twoside]{article}
\usepackage[latin1]{inputenc}
\usepackage[T1]{fontenc}
\usepackage{parskip}
\usepackage{amsmath,amsfonts,amssymb,amsxtra}
\usepackage{latexsym}
\usepackage{theorem}
\usepackage{graphicx,psfrag}
\usepackage{verbatim}
\usepackage{hyperref}

\newtheorem{theo}{Theorem}
\newtheorem{prop}{Proposition}
\newtheorem{lemma}{Lemma}
\newtheorem{cor}{Corollary}

{\theoremstyle{plain}                       
\theorembodyfont{\rmfamily}
\newtheorem{Remark}{Remark}}

{\theoremstyle{plain}                       
\theorembodyfont{\rmfamily}
\newtheorem{Example}{Example}}

{\theoremstyle{plain}                       
\theorembodyfont{\rmfamily}
\newtheorem{Definition}{Definition}}

{\theoremstyle{plain}                       
\theorembodyfont{\rmfamily}
\newtheorem{Condition}{Condition}}

\def\M{\mathcal{M}}
\def\F{\mathcal{F}}
\def\B{\mathcal{B}}
\def\I{\mathcal{I}}
\def\S{\Sigma}
\def\Sp{\Sigma^+}
\def\BS{\B(\S)}
\def\N{\mathbb{N}}
\def\R{\mathbb{R}}
\def\Z{\mathbb{Z}}

\begin{document}

\title{Contractive Markov systems II}

\author{Ivan Werner\\
{\small Email: ivan\_werner@mail.ru}} \maketitle

\begin{abstract}\noindent
Discrete time random dynamical systems with countably many maps which admit countable Markov partitions on complete metric spaces such that the resulting Markov systems are uniform continuous and contractive are considered.  A notion of a generating communication class of such a system is introduced, which includes every communication class if the system has a finite Markov partition. It is shown that the ergodic decomposition of an equilibrium state associated with such a system is purely atomic and can be exhaustively described using the generating communication classes if the system satisfies  an absolute continuity condition (ACC).   In such a case, each invariant Borel probability measure which is an image of an ergodic component of an equilibrium state under the coding map can be obtained by a random walk starting at any point in the corresponding generating communication class. As a by-product, a practical method for a computation of the entropy of the equilibrium states is obtained. Finally, it is shown  that such a non-degenerate system  satisfying  the ACC which in addition has a dominating Markov chain and a finite \eqref{cgc} has a unique invariant Borel probability measure if and only if it has a single generating communication class. Some sufficient conditions for the ACC are provided. 


\end{abstract}

\tableofcontents

\section{Introduction}

The purpose of this article is, in particular, to provide criteria for the uniqueness of an invariant Borel probability measure  for the random dynamical systems introduced in  \cite{Wer1} as {\it contractive Markov systems} (the special case of the systems when all its  maps are contractions on a compact metric space was first considered in \cite{Wer0}). They consistently extend and  unify several previously studied structures, such as  {\it weighted directed graphs} (discrete  homogeneous {\it Markov chains} \cite{M}), {\it random systems with compete connections} \cite{OM}, \cite{DF}, {\it learning models} \cite{Ka}, \cite{Is}, {\it $g$-functions} \cite{ce} and {\it  iterated function systems  with place dependent probabilities (IFSPDP)} \cite{BDEG}, \cite{Elton}. Such a generalization allows to extend the  powerful and practical language of weighted directed graphs  based on the notion of a communication class and combine it with the general logic of ergodic theory and thermodynamic formalism  to formulate new insights on the behaviour of such random dynamical systems. It is not surprising that it also turns out to be  a source of challenging and enlightening examples for other areas of mathematics, not excluding ergodic theory and thermodynamic formalism, as it is demonstrated in \cite{Wer6}, \cite{Wer10}, \cite{Wer12} and in this article, because the same role has already been  played by some of its particular cases such as weighted directed graphs and $g$-functions. Not to ignore is also the applicative power inherited by such systems from weighted directed graphs and IFSPDP, which has been enjoyed and challenged by many scientific applications (e.g. see one of the latest which has challenged IFSPDP  in \cite{FHK}). It is drastically magnified by a much richer spectrum of processes which can be generated by contractive Markov systems, compared to the trivial case of weighted directed graphs, which is due to a much greater complexity of the topological structure of such systems, whereas the algorithm for the generation of the processes remains the same. Probably, because of this complexity, the behaviour of such systems is still not completely clarified.

The uniqueness of a stationary state for such a finite irreducible system was proved in  \cite{Wer1} under the condition that the partition of the system on a locally compact metric space consists of open sets and the probability functions are Dini-continuous (have a summable variation) and bounded away from zero. This already covered the corresponding results for finite weighted directed graphs and $g$-functions \cite{Le}, \cite{Wal}, but its applicability still remained restricted to disconnected state spaces. The proof of the result in  \cite{Wer1} was an extension of the proof which had been given by M. F. Barnsley et al. in \cite{BDEG} and \cite{BDEGE} for IFSPDP, which had used the classical coupling method. The result then was extended  by K. Horbacz and T. Szarek on Polish spaces \cite{HS} by applying some previous results obtained by the second author for Markov operators on Polish spaces \cite{Sz}.  The method which was chosen in \cite{HS} also required that the Markov operator associated with the system mapped continuous functions on continuous functions, and this forced the authors to keep the assumption of the openness on the Markov partition.

The main obstacle which is associated with an arbitrary Markov partition lies in the proof of the existence of an invariant measure.
Recently, it has been overcome by identifying the conditions for the existence of  {\it equilibrium states} or, in general, {\it asymptotic states} on the code space associated with such a system,  which are then mapped on  invariant measures by a coding map \cite{Wer11}. (The development of the approach began with the construction of the coding map in \cite{Wer3} and the connection of the invariant measures with the equilibrium states for a {\it local energy function}, obtained by means of the coding map, in \cite{Wer6}.) In spite of the fact that the problem turned out to be beyond the current  theory of equilibrium states, even in the case of finite systems with  open Markov partitions \cite{Wer6}, the existence of the invariant Borel probability measures has been proved in \cite{Wer11} for such systems with some proper Markov partitions on complete metric spaces  satisfying a {\it non-degeneracy} or a {\it consistency} condition, which successively weaken  the openness condition. In particular, the consistency condition is satisfied by all random dynamical systems with continuous maps and probability functions which have  finitely many uniformly continuous restrictions on each atom of their Markov partitions  (the case with infinitely many restrictions of the probability functions on some atoms requires, in addition, the existence of a {\it dominating  Markov chain} \eqref{dmcc}, see  \cite{Wer11}). 

In the case of a countably infinite Markov partition, the existence of an equilibrium state on the code space requires an additional condition (Condition 2  in \cite{Wer11}), which corresponds to the {\it positive recurrence} in the case of a discrete homogeneous Markov chain. It is automatically satisfied in a finite case. It was shown in \cite{Wer11} that the condition is necessary and sufficient for the existence of an invariant measure for a uniformly continuous,  contractive, Markov system in the non-degenerate case and sufficient in the consistent case provided the system satisfies  \eqref{cgc} and \eqref{dmcc}. This condition leads to the definition of a {\it generating communication class} in this article (Section \ref{gccs}), and allows us to associate the ergodic components of an equilibrium state with the generating communication classes supporting the image of the equilibrium state under the coding map (Section \ref{deds}) if the system satisfies an absolute continuity condition  (Condition \ref{acC}) introduced in Section \ref{fcs}. In particular, it is shown that such a uniformly continuous, contractive, Markov system satisfying  Conditions \ref{acC} and \ref{dtc} has a unique equilibrium state in the non-degenerate case if and only if it has a single generating communication class. This automatically translates into the necessary and sufficient condition for the uniqueness of an invariant Borel probability measure for such a system   in the non-degenerate case (Section \ref{cim}), through the one-to-one correspondence of the invariant measures and equilibrium states established in \cite{Wer11}. Moreover, in the case with several generating communication classes, it is shown that each such class supports a unique invariant ergodic measure,  the values of which can be computed through an ergodic average of a random walk started at any point in that class (Theorem \ref{esgc}).  As an application,  a practical method for a computation of Shannon-Kolmogorov-Sinai entropy of processes generated by such systems on the code space is obtained (Corollary \ref{Ec}).

The key to the proof is the ergodic decomposition of equilibrium states associated with the considered random dynamical systems. It turns out that the well-known result from the theory of equilibrium states is not applicable to the systems in the infinite case (see Section \ref{edes} for an explanation). A self-contained proof of the ergodic decomposition of equilibrium states for such systems is given in  Section \ref{edes}.

Finally, Section \ref{scac} provides some sufficient conditions for  Condition \ref{acC}, which do not require the boundedness away from zero of the probability functions.

It was pointed out by an anonymous reviewer that it might be appropriate here to cite the work of Ch. Walkden \cite{Wak}, where, in particular,  the stability of a unique invariant probability measure with respect to a change of the probability functions is studied, and the works of O. Sarig \cite{Sa1}, \cite{Sa2} and \cite{Sa3}, where, in particular, a notion of a positive recurrence for potentials on countable Markov shifts is introduced and the thermodynamic formalism for the potentials with a summable variation is developed (note that the potentials on countable Markov shifts associated with the random dynamical systems considered in this article are not even upper semicontinuous in general, even in the case with an open Markov partition).

\section{Definitions and notation}

In this article, a {\it random dynamical system} on a metric space $(K, d)$ is a family $\mathcal{D}_R:=(K, w_e,p_e)_{e\in E'}$ where $E'$ is an at most countable set,  $(w_e)_{e\in E'}$ is a family of Borel-measurable  {\it maps} of $K$ into itself and $(p_e)_{e\in E'}$ is a family of Borel-measurable {\it probability functions} $p_e:K\longrightarrow [0,1]$.  The random dynamical system acts on the set of all real-valued non-negative Borel-measurable function on $K$, $\mathcal{L}(K)$,  by a {\it Markov operator} $U$ given by
\[Uf:=\sum\limits_{e\in E'}p_ef\circ w_e\]
for all $f\in\mathcal{L}(K)$ and on the set of  Borel probability measures on $K$, $P(K)$, by the {\it adjoint operator} $U^*$ given by
$U^*\nu(f):=\int Ufd\nu$ for all  bounded $f\in \mathcal{L}(K)$ and $\nu\in P(K)$. We say that $\mu\in P(K)$ is an {\it invariant measure} of $\mathcal{D}_R$ if and only if $U^*\mu = \mu$. 

It is clear from the definition of $U$ that each $w_e$ needs to be defined only on the set $\{p_e>0\}$ (in such a case, $w_e$ can be considered to be extended on $K$ arbitrarily).   A  random dynamical system $(K, w_e,p_e)_{e\in E}$ is called a {\it Markov system} if and only if there exists a partition of $K$ into non-empty Borel subsets $(K_i)_{i\in N}$ (the case where $N$ has only one element is not excluded)  such that for every $e\in E$ there exist $i(e),t(e)\in N$ such that $\emptyset\neq\{p_e>0\}\subset K_{i(e)}$  and $w_e(K_{i(e)})\subset K_{t(e)}$. $i:E\longrightarrow N$ is required to be surjective.  Clearly, this defines a topological structure on $K$ which generalizes a weighted directed graph. $\mathcal{D}_R$ is said to have a {\it Markov partition} if and only if the restrictions of its probability functions and maps on the atoms of a partition form a Markov system (after a possible enlargement of the index set $E'$) . The atoms of the partition are called the {\it vertex sets} of the Markov system, and $N$ is called the {\it set of vertices}.  We will denote a Markov system by $\M:=(K_{i(e)},w_e,p_e)_{e\in E}$.  $\M$ is called {\it proper} if and only if $N$ has more than one element, and {\it countable} if and only if $N$ is at most countable. We call $\M$ {\it positive} if and only if $p_e|_{K_{i(e)}}>0$ fot all $e\in E$. $\M$ is called {\it (uniformly) continuous} if and only if restrictions of the maps $w_e|_{K_{i(e)}}$ and the probability  functions  $p_e|_{K_{i(e)}}$ are (uniformly) continuous.  $\M$ is called {\it contractive} if and only if there exists $0<a<1$ such that
 \begin{equation}\label{cac}
 \sum\limits_{e\in E, i(e)=j}p_e(x)d(w_ex,w_ey)\leq ad(x,y)\mbox{ for all }x,y\in K_j\mbox{ and } j\in N.
 \end{equation}
The condition goes back to R. Isaac \cite{Is} for the case when $N$ has one element.

If $\M$ is uniformly continuous, let $\bar p_e$ denote the continuous extension of $p_e|_{K_{i(e)}}$ on the closure $\bar K_{i(e)}$ which is extended on $K$ by zero, and $\bar w_e$ denote the continuous extension of $w_e|_{K_{i(e)}}$ on $\bar K_{i(e)}$ which is extended on $K$ arbitrarily.

A sequence $(e_1,e_2,...,e_n)\in E^n$ is called a {\it path} of the Markov systems if and only if $i(e_{i+1}) = t(e_i)$ for all $1\leq i\leq n-1$. In such a case, $n$ is called the {\it length of the path}.
 We say that $j\in N$ is {\it accessible} from $i\in N$ if and only if either $i=j$ or there  exists a path $(e_1,...,e_n)$ such that $i = i(e_1)$ and $j = t(e_n)$. We say that two  vertices {\it communicate} if and only if one is accessible from the other and vice versa.  $i\in N$ is called {\it essential} if and only if $i$ communicates with every $j\in N$ which is accessible from $i$. Clearly communication is an equivalence relation on $N$, and therefore $N$ splits into equivalence classes. Let $\bigcup_{i\in I}c_i\subset N$ be the partition of the set of all essential vertices of $\M$ into equivalence classes. For each $i\in I$, set
\[C_i:=\bigcup\limits_{j\in c_i}K_j.\]
We will call both $c_i$ and $C_i$ an  {\it(irreducible) communication class} of $\M$  if no confusion is possible. 

Let $E$, $N$ and $I$  be provided with the discrete topologies. Let $\bar E:=E\cup\{\infty\}$ denote the one-point compactification of $E$, and set $w_\infty:= id$, $p_\infty:=0$. Fix $n_\infty\in N$, and set $i(\infty)=t(\infty)=n_\infty$. Let
$\Sigma:=\{(...,\sigma_{-1},\sigma_0,\sigma_1,...):\sigma_i\in \bar E\
\forall i\in\mathbb{Z}\}$ and $\Sigma^+:=\{(\sigma_1,\sigma_2,...):\sigma_i\in \bar E\
\forall i\in\mathbb{N}\}$ provided with the product topologies and Borel $\sigma$-algebras.  $_m[e_m,...,e_n]:=\{\sigma\in\Sigma|\ \sigma_i=e_i\mbox{ for all }m\leq i\leq n\}$ is called a {\it cylinder set}. A cylinder set in $\Sigma^+$ will be denote by $_k[e_k,...,e_n]^+$.
For $m\leq 0$, let $\F_m$ denote the  $\sigma$-algebra on $\Sigma$ generated by the cylinder sets of the form $_m[e_m,...,e_0]$, and $\F$ be the $\sigma$-algebra generated by $\bigcup_{m\leq 0}\F_m$.
Let $S:\Sigma\longrightarrow\Sigma$ denote the left {\it shift map}, given by $(S\sigma)_{i-1}=\sigma_i$ for all $i\in\mathbb{Z}$.
Set
\[\Sigma_G :=\left\{\sigma\in\Sigma |\  i(\sigma_{n+1}) = t(\sigma_n),\ \sigma_n\in E\mbox{ for all }n\in\mathbb{Z}\right\},\]
and analogously $\Sigma^+_G$. $\Sigma_G$ is called the {\it path space} of $\M$. 

We will denote the Borel $\sigma$-algebra of on a topological space $X$ by $\B(X)$ and the set of all Borel probability measures on $X$ by $P(X)$. Let $P(\M)\subset P(K)$ denote the set all invariant measures of $\M$, and $P_S(\S)\subset P(\S)$ the set of all shift-invariant measures.

 For $x\in K$, let $P_x\in P(\Sigma^+)$ be given by
\[P_x\left( _1[e_1,...,e_k]^+\right):=p_{e_1}(x)p_{e_2}(w_{e_1}x)...p_{e_k}(w_{e_{k-1}}\circ...\circ w_{e_1}x)\]
for all $_1[e_1,...,e_k]^+\subset\Sigma^+$. We call $P_x$ a {\it path measure} of $\M$. For any $_m[e_m,...,e_k]\subset\Sigma$, set
$P^m_x(_m[e_m,...,e_k]):=P_x(_1[e_m,...,e_k]^+)$. For $\mu\in P(\M)$, $\Phi(\mu)\in P_S(\Sigma)$ is uniquely defined by 
\[\Phi(\mu)\left(_m[e_m,...,e_k]\right):=\int P^m_x(_m[e_m,...,e_k])d\mu(x)\]
for all $_m[e_m,...,e_k]\subset\Sigma$.
 For $\nu\in P(K)$,  $\tilde \phi(\nu)$ denotes the probability measure on  the product $\sigma$-algebra $\B(K)\otimes
\B(\Sigma^+)$ given by
\[\tilde \phi(\nu)(\Omega)=\int P_x\left(\left\{\sigma\in\Sigma^+:(x,\sigma)\in\Omega\right\}\right)d\nu(x)\]
for all $\Omega\in \B(K)\otimes
\B(\Sigma^+)$.

For a measurable map on a measure space $f:(X, \mathcal{A}, \lambda)\longrightarrow(Y,\B)$, let $f(\lambda)$ denote the measure on $\B$ given by $f(\lambda)(B):=\lambda(f^{-1}(B))$ for all $B\in \B$. As usual, let $\Delta$ denote the symmetric difference for sets and $\ll$ denote the absolute continuity relation for measures.

For $M\in P_S(\S)$, let 
\[\mathcal{I}_M:=\left\{B\in\B\left(\S\right)|\ M\left(S^{-1}B\Delta B\right)=0\right\}.\]
It is easy to check that $\mathcal{I}_M$ is a sub-$\sigma$-algebra.

\section{Results}

Let  $\M:=(K_{i(e)},w_e, p_e)_{e\in E}$ be a countable Markov system on a complete metric space $(K,d)$.
For each $i\in N$, fix $x_i\in K_i$, and set
\[D:=\left\{\sigma\in\Sigma_G\left|\ \lim\limits_{m\to-\infty}w_{\sigma_0}\circ ... \circ w_{\sigma_m}(x_{i(\sigma_m)})\mbox{ exists}\right.\right\}.\]
For every $\sigma\in\Sigma$, set
\begin{equation*}
    F(\sigma):=\left\{\begin{array}{cc}
    \lim\limits_{m\to-\infty}w_{\sigma_0}\circ w_{\sigma_{-1}}\circ...\circ w_{\sigma_{m}}(x_{i(\sigma_{m})})&  \mbox{if }\sigma\in D\\
     x_{t(\sigma_0)}& \mbox{ otherwise, }
     \end{array}\right..
     \end{equation*}
$F$ is called the {\it coding map} of the Markov system. Clearly, it is  $\F$-Borel-measurable.   Set
 \[E(\M):=\left\{\Lambda\in P_S(\Sigma)|\  \Lambda(D) = 1\mbox{ and }E_\Lambda(1_{_1[e]}|\F)=p_e\circ F\  \Lambda\mbox{-a.e. for all }e\in E\right\}.\]
The following non-degeneracy condition was introduced in \cite{Wer11}. It admits some proper Markov systems on connected spaces. 
\begin{Definition}
     Set $T_j:=\{\sigma\in\Sigma_G|\ t(\sigma_0)=j\}$ for $j\in N$.  Suppose $p_e|_{K_{i(e)}}$ is uniformly continuous for all $e\in E$. Let $\tilde E(\M)$ denote the set 
\[\left\{\Lambda\in P_S(\Sigma)|\  \Lambda(D) = 1\mbox{ and }E_\Lambda(1_{_1[e]}|\F)=\bar p_e\circ F1_{T_{i(e)}}\  \Lambda\mbox{-a.e. for all }e\in E\right\}.\]
We call  $\M$ {\it non-degenerate} if and only if  for every $\Lambda\in \tilde E(\M)$ there exists $i\in N$ such that $\Lambda(T_i\cap F^{-1}(K_i))>0$. Clearly, every uniformly continuous Markov system with an open partition is non-degenerate, as always $T_i\subset F^{-1}(\bar K_i)$ for all $i\in N$ (see Example \ref{se} and \cite{Wer11} for some examples of proper non-degenerate Markov systems  on connected spaces). By Theorem 2 in \cite{Wer11}, the non-degeneracy  is equivalent to  $\tilde E(\M)=E(\M)$. Also, in \cite{Wer11}, a practical sufficient condition for the non-degeneracy is given.
\end{Definition}

\subsection{The ergodic decomposition of equilibrium states}\label{edes}

It is a well-known fact from the theory of equilibrium states that an equilibrium state  $M$ of a continuous dynamical system  on a compact metric space $X$ with the upper semi-continuous entropy function (and therefore, with a finite topological entropy), minimising the free energy, for an upper semi-continuous energy function $\psi:X\longrightarrow[-\infty,+\infty)$ decomposes into ergodic components such that $M$-almost every of them is again an equilibrium state for $\psi$, e.g. see Theorem 4.3.9 in \cite{Kel}. 

If $E$ is finite, then, by Theorem  1 in \cite{Wer11}, every $M\in E(\M)$ is an equilibrium state in the thermodynamic sense for an energy function $u:\S\longrightarrow [-\infty,0]$ given by
 \begin{equation*}
    u(\sigma):=\left\{\begin{array}{cc}
    \log  p_{\sigma_1}\circ F(\sigma)&  \mbox{if }\sigma\in D\\
    -\infty& \mbox{ otherwise }
     \end{array}\right.\ \ \ \mbox{ for all } \sigma\in\Sigma
     \end{equation*}
   (with the definition $\log(0) := -\infty$), and therefore, Theorem 4.3.9 in \cite{Kel} can be applied to it (though, $u$ is not upper semi-continuous, but $u\in\mathcal{L}^1(M)$ and the proof of Theorem 4.3.9 in \cite{Kel} applies to it in this case also).

In the case of a countably infinite $E$, the definition of an equilibrium state from \cite{Kel} does not extend to $u$, as $u$ is not bounded from below and the Kolmogorov-Sinai entropy of an invariant Borel probability measure on $\S$ can be infinite. The definition of equilibrium states in the thermodynamic sense for $u$ which was given in \cite{Wer11} covers, by Theorem  1 in \cite{Wer11}, only the members of $E(\M)$ with finite entropy, but an example where all members of $E(\M)$ have infinite entropy was given in \cite{Wer11} also. The author is not aware of any result on the ergodic decomposition of $M\in E(\M)$ which could be applied in this case. Such a result is provided in this subsection.

We will use the following well-known theorem on the existence of the ergodic decomposition of invariant measures. 
\begin{theo}\label{ted}
    Let $M\in P_S(\S)$. Then for every $\omega\in\S$ there exists  $\Lambda_\omega\in P_S(\S)$ such that the following holds true.\\
(i) For every $f\in\mathcal{L}^1(M)$, $\omega\mapsto\int f\ d\Lambda_\omega$ is $\I_M$-measurable and $\int f\ d\Lambda_\omega = E_M(f|\mathcal{I}_M)(\omega)$ for $M$-a.e. $\omega$. \\
 (ii) $\Lambda_{S\omega}=\Lambda_\omega$ for all $\omega\in\S$, and\\
(iii)  $\Lambda_\omega$ is ergodic for all $\omega\in\S$.
\end{theo}
{\it Proof.} 
The assertion follows by  Theorem 2.3.3 in \cite{Kel}.
\hfill$\Box$

\begin{Definition}
The collection of measures $\{\Lambda_\omega\}_{\omega\in\S}$ with the properties from Theorem \ref{ted} is called {\it the ergodic decomposition of} $M$, since it is unique $M$-a.s., and is denoted by $M=\int\Lambda_\omega dM(\omega)$.
\end{Definition}

The following lemma is  well known, unfortunately, the author didn't find anything to cite. 
\begin{lemma}\label{psa}
   Let $M\in P_S(\S)$ and $Q\in \I_M$. Then there exists $A\in\F$ such that $M(Q\Delta A)=0$.
\end{lemma}
{\it Proof.} 
The proof is a straightforward exercise.
\hfill$\Box$

The following theorem is probably also well-known, but the author didn't find anything to cite.
\begin{theo}\label{cmedt}
    Let $M\in P_S(\S)$ and $M=\int \Lambda_\omega dM(\omega)$ be its ergodic decomposition. Let $e\in E$ and $f_e$ be a  version of $E_M(1_{_1[e]}|\F)$. Then, for $M$-a.e. $\omega\in\S$,
\[E_{\Lambda_\omega}(1_{_1[e]}|\F)=f_e\ \ \ \Lambda_\omega\mbox{-a.e. }.\]
\end{theo}
{\it Proof.} 
Let $A\in\F$ and $Q\in\I_M$. Then, by Lemma \ref{psa},  there exists $\tilde Q\in\F$ such that $M(Q\Delta\tilde Q)=0$.  Therefore, since $E_M(1_{_1[e]}|\F)=f_e$ $M$-a.e.,
\begin{eqnarray*}
  &&\int\limits_{Q}\int\limits_A1_{_1[e]}d\Lambda_\omega dM(\omega)=\int\limits_{Q}1_A1_{_1[e]}dM=\int\limits_{\tilde Q}1_A1_{_1[e]}dM=\int\limits_{\tilde Q}1_Af_edM\\
&=&\int\limits_{Q}1_Af_edM=\int\limits_{Q}\int\limits_Af_ed\Lambda_\omega dM(\omega).
\end{eqnarray*}
Hence, since $Q\in\I_M$ was arbitrary,  and $\omega\mapsto\int .\Lambda_\omega$ is $\I_M$-measurable,
\[\int\limits_A1_{_1[e]}d\Lambda_\omega=\int\limits_Af_ed\Lambda_\omega\ \ \ \mbox{ for }M\mbox{-a.e. }\omega\in\S.\]
Let $\mathcal{G}\subset\F$ denote the collection of cylinder sets of the form $_m[e_m,...,e_0]$, $e_m,...,e_0\in\bar E$ and $m\leq 0$. Since $\mathcal{G}$ is countable, there exists $X\in\BS$ with $M(X)=1$ such that for every $\omega\in X$,
 \begin{equation}\label{esce}
    \int\limits_A1_{_1[e]}d\Lambda_\omega=\int\limits_Af_ed\Lambda_\omega\mbox{ for all }A\in\mathcal{G},
 \end{equation}
i.e. for every $\omega\in X$ the measures on $\F$ given by the left hand  and the right hand sides of \eqref{esce}  agree on $\mathcal{G}$. Since $\mathcal{G}$ generates $\F$, is $\cap$-stable and covers $\S$, they agree also on $\F$. Thus, for every $\omega\in X$,
\[E_{\Lambda_\omega}(1_{_1[e]}|\F)=f_e\ \ \ \Lambda_\omega\mbox{-a.e. }.\]
\hfill$\Box$

\begin{cor}\label{esdc}
    Suppose $\M$ is uniformly continuous. Let $M\in P_S(\S)$ and $M=\int \Lambda_\omega dM(\omega)$ be its ergodic decomposition.\\ (i) If $M\in E(\M)$, then $\Lambda_\omega\in  E(\M)$ for $M$-a.e. $\omega\in\S$.\\
(ii) If $M\in \tilde E(\M)$, then $\Lambda_\omega\in\tilde  E(\M)$ for $M$-a.e. $\omega\in\S$.
\end{cor}
{\it Proof.} 
(i)  Let $e\in E$. Since $E_M(1_{_1[e]}|\F)=p_e\circ F$ $M$-a.e., by Theorem \ref{cmedt}, there exists $X\in\BS$ with $M(X)=1$ such that  for every $\omega\in X$,
\[E_{\Lambda_\omega}(1_{_1[e]}|\F)=p_e\circ F\ \ \ \Lambda_\omega\mbox{-a.e. }.\]
Since $1=M(D)=\int\Lambda_\omega(D)dM(\omega)$, there exists $Y\in\BS$ with $M(Y)=1$ such that $\Lambda_\omega(D)=1$ for all $\omega\in Y$. Thus $\Lambda_\omega\in E(\M)$ for all $\omega\in X\cap Y$.

The proof of (ii) is the same. One  needs only  to replace $p_e\circ F$ with $\bar p_e\circ F1_{T_{i(e)}}$.
\hfill$\Box$

\subsection{Generating points and communication classes}\label{gccs}

\begin{Definition}\label{gpd}
Let $N$ be provided with  Borel $\sigma$-algebra. We call $x_0\in K$ a {\it generating point } (for $\M$)  if and only if the sequence of probability measures $(\alpha^{x_0}_n)_{n\in\mathbb{N}}$ on $N$ given by 
\begin{equation}\label{gpm}
\alpha^{x_0}_n\{j\}:=\frac{1}{n}\sum\limits_{k=1}^n{U^*}^k\delta_{x_0}\left(K_j\right)\ \ \ \mbox{ for all }j\in N\mbox{ and }n\in\mathbb{N}
\end{equation}
is uniformly tight. We call $i_0\in N$ {\it generating} (for $\M$) if and only if there exists a generating point $x_0\in K_{i_0}$. Clearly, every $x\in K$ is generating if $N$ is finite.
\end{Definition}

\begin{lemma}\label{gpil}
    (i) Suppose $\M$ is positive, $i_0\in N$ is {\it generating} for $\M$, and  $j_0\in N$ is accessible from $i_0$. Then $j_0$ is also generating.\\
(ii) Suppose $\M$ is positive. Let $i\in N$ be essential and $C_i\subset N$ be the communication class containing $i$. If $i$ is generating, then every $j\in C_i$ is generating.
\end{lemma}
{\it Proof.} 
(i) Let $x_0\in K_{i_0}$ be generating. Let $(e'_1,...,e'_m)$ be a path such that $i_0=i(e'_1)$ and $j_0=t(e'_m)$. Set $y:=w_{e'_m}\circ...\circ w_{e'_1}(x_0)$. Then $y\in K_{j_0}$.  Let $\epsilon>0$.
Set $\delta:=P_{x_0}(\L_1[e'_1,...,e'_m]^+)$. Then, by the hypothesis, $\delta>0$, and there exists finite $V_\epsilon\subset N$ such that $\alpha^{x_0}_n(N\setminus V_\epsilon)<\epsilon\delta/2$ for all $n\in\N$.  Choose $n_0>m$ such that $m/n_0<\epsilon/2$. Observe that for all $n>m$ and $j\in N$,
\begin{eqnarray*}
\alpha^{x_0}_n\{j\}&=&\frac{1}{n}\sum\limits_{k=1}^n\int 1_{K_j}\circ w_{\sigma_k}\circ...\circ w_{\sigma_1}(x_0)dP_{x_0}(\sigma)\\
&\geq&\frac{1}{n}\sum\limits_{k=m+1}^n\int\limits_{_1[e'_1,...,e'_m]} 1_{K_j}\circ w_{\sigma_k}\circ...\circ w_{\sigma_{m+1}}(y)dP_{x_0}(\sigma)\\
&=&\frac{1}{n}\sum\limits_{k=m+1}^n\sum\limits_{e_{m+1},...,e_k}P_{x_0}\left( _1[e'_1,...,e'_m]^+\right)P_{y}\left( _1[e_{m+1},...,e_k]^+\right)\\
&&\ \ \ \ \ \ \ \ \ \ \ \ \ \ \ \ \ \ \ \ \ \times1_{K_j}\circ w_{e_k}\circ...\circ w_{e_{m+1}}(y)\\
&=&\delta\frac{1}{n}\sum\limits_{k=m+1}^n\int1_{K_j}\circ w_{\sigma_{k-m}}\circ...\circ w_{\sigma_{1}}(y)P_{y}(\sigma)\\
&=&\delta\alpha^y_n\{j\}-\delta\frac{1}{n}\sum\limits_{k=n-m+1}^{n}\int1_{K_j}\circ w_{\sigma_{k}}\circ...\circ w_{\sigma_{1}}(y)P_{y}(\sigma).
\end{eqnarray*}
Therefore, for all $n\geq n_0$,
\begin{eqnarray*}
\alpha^y_n(N\setminus V_\epsilon)&\leq\frac{1}{\delta}\alpha^{x_0}_n(N\setminus V_\epsilon)+\frac{m}{n}<\frac{\epsilon}{2}+\frac{\epsilon}{2}=\epsilon.
\end{eqnarray*}
For each $1\leq n<n_0$, choose $V_n\subset N$ such that $\alpha^y_n(N\setminus V_n)<\epsilon$, and set $V'_\epsilon:=V_\epsilon\cup\bigcup_{n=1}^{n_0-1}V_n$. Then $V'_\epsilon$ is finite, and $\alpha^y_n(N\setminus V'_\epsilon)<\epsilon$ for all $n\in\N$. This completes the proof of (i).

(ii) follows immediately from (i). 
\hfill$\Box$

By Lemma \ref{gpil} (ii), we can make the following definition.
\begin{Definition}
Suppose $\M$ is positive. Let $i\in I$. We call $c_i$ and $C_i$ {\it generating} if and only if there exists $j\in c_i$ which is generating.
Let $I_G\subset I$ denote the set of all generating communication classes of $\M$.
\end{Definition}

\subsection{The absolute continuity condition (ACC)}\label{fcs}

The following condition for a Markov system seems to be of a fundamental nature. Apparently, it was first shown to be satisfied, for the case of a single vertex set, in a work by J. Elton \cite{Elton}  in the case where each $p_e$ is Dini-continuous (has a summble variation) and  bounded away from zero.
\begin{Condition}[ACC]\label{acC}
 For all $i\in N$ and $x,y\in K_i$, 
\begin{equation*}
P_x\ll P_y.
\end{equation*}
\end{Condition}
Note that  Condition \ref{acC} implies, in particular, that $\M$ is positive. The following theorem collects some well-known equivalent statements to Condition \ref{acC}  if $\M$ is positive, which will be used later. For some sufficient conditions for Condition \ref{acC}, see  Section \ref{scac}.

\begin{Definition}
Suppose $\M$ is positive. Let $i\in N$ and $x,y\in K_i$. For $k\in\N$, let $\B_k$ denote the $\sigma$-algebra on $\Sp$ generated by cylinder sets of the form $_1[e_1,...,e_k]^+$. Let $P_x|_{\B_k}$ denote the restriction of $P_x$ on $\B_k$. Since $M$ is positive, there exist Radon-Nikodym  derivatives \[Z^{xy}_k:=\frac{dP_x|_{\B_k}}{dP_{y}|_{\B_k}}\] for all $k\in\N$. It is easy to check that $(Z^{xy}_k, \B_k)_{k\in\N}$ is a $P_{y}$-martingale with $\int Z^{xy}_k dP_{y} = 1$ for all $k\in\N$. Hence, by Doob's Martingale Theorem, $Z^{xy}_\infty:=\lim_{k\to\infty}Z^{xy}_k$ exists $P_{y}$-a.e. and is integrable. 
\end{Definition}

\begin{theo}\label{acet}
    Suppose $\M$ is positive. Let $i\in N$ and $x,y\in K_i$. Then the following are equivalent.\\
(i) $P_x\ll P_y$. \\
(ii) $\int Z^{xy}_\infty dP_y= 1$.\\
(iii) $P_x\{ Z^{xy}_\infty<\infty\} =1$.\\
 (iv) $(Z^{xy}_k)_{k\in\N}$ is  uniformly integrable with respect to $P_y$.\\
 (v) $\sup_{k\in\N}P_x\{Z^{xy}_k>c\}\to 0$ as $c\to\infty$.\\
 (vi) $\sup_{k\in\N}P_x\{\log Z^{xy}_k>c\}\to 0$ as $c\to\infty$.
\end{theo}
{\it Proof.} 
  The equivalence of (i), (ii) and (iii) follows  by Shiryaev's Local Absolute Continuity Theorem, e.g. Theorem 2, p.514, in \cite{Shi}.

Assertion (ii) means  $\int Z^{xy}_\infty dP_{y}=\lim_{k\to\infty}\int Z^{xy}_k dP_{y}$, and since  $ Z^{xy}_k\to   Z^{xy}_\infty$ $P_y$-a.e., the convergence of the integrals is equivalent to $(Z^{xy}_k)_{k\in\N}$ being uniformly integrable with respect to $P_y$, e.g. Theorem 5 p.205 in \cite{Shi}. 

The equivalence of (v) and (iv) follows from the definition of $Z^{xy}_k$, since 
\[\int\limits_{\{Z^{xy}_k>c\}}Z^{xy}_kdP_y=P_x\{Z^{xy}_k>c\}\]
for all $k\in\N$ and $c>0$.

The equivalence of (v) and (vi) is obvious.
 \hfill$\Box$

See \cite{KLS} for further necessary and sufficient conditions for the absolute continuity of measures and Section 4.5, Vol. 1 in \cite{Bog} for further necessary and sufficient conditions for the uniform integrability.

\begin{lemma}\label{tgl}
    Suppose $\M$ satisfies Condition \ref{acC}. Let $i_0\in N$,  and $x_0\in K_{i_0}$ be {\it generating}. Then  every $x\in K_{i_0}$  is generating.
\end{lemma}
{\it Proof.} 
Let $x\in K_{i_0}$ and  $n\in\N$. By Theorem \ref{acet}, $P_x\ll P_{x_0}$ is equivalent to the uniform integrability of $(Z^{xx_0}_k)_{k\in\N}$ with respect to $P_{x_0}$.  Let $\epsilon>0$. Choose $c>0$ such that
\[\sup\limits_{k\in\N}\int\limits_{\left\{Z^{xx_0}_k>c\right\}}Z^{xx_0}_kdP_{x_0}<\frac{\epsilon}{2}.\]
Choose finite $V_\epsilon\subset N$ such that $\alpha^{x_0}_n(N\setminus V_\epsilon)<\epsilon/(2c)$ for all $n\in\N$. Observe that for every $j\in\N$ and $k\in\N$,
\begin{eqnarray*}
{U^*}^k\delta_x(K_j)&=&\int 1_{K_j}\circ w_{\sigma_k}\circ...\circ w_{\sigma_1}(x)dP_{x}(\sigma)\\
&=&\int 1_{K_j}\circ w_{\sigma_k}\circ...\circ w_{\sigma_1}(x_0)dP_{x}(\sigma)\\
&=&\int 1_{K_j}\circ w_{\sigma_k}\circ...\circ w_{\sigma_1}(x_0)Z^{xx_0}_k(\sigma)dP_{x_0}(\sigma)\\
&\leq&c\int 1_{K_j}\circ w_{\sigma_k}\circ...\circ w_{\sigma_1}(x_0)dP_{x_0}(\sigma)\\
&&+\int\limits_{\left\{Z^{xx_0}_k>c\right\}} 1_{K_j}\circ w_{\sigma_k}\circ...\circ w_{\sigma_1}(x_0)Z^{xx_0}_k(\sigma)dP_{x_0}(\sigma).
\end{eqnarray*}
Hence,
\begin{eqnarray*}
\alpha^x_n(N\setminus V_\epsilon)&\leq&c\alpha^{x_0}_n(N\setminus V_\epsilon) +\frac{1}{n}\sum\limits_{k=1}^n\int\limits_{\left\{Z^{xx_0}_k>c\right\}} Z^{xx_0}_k(\sigma)dP_{x_0}(\sigma)\\
&\leq&\frac{\epsilon}{2} +\sup\limits_{k\in\N}\int\limits_{\left\{Z^{xx_0}_k>c\right\}} Z^{xx_0}_k(\sigma)dP_{x_0}(\sigma)\\
&<&\epsilon.
\end{eqnarray*}
\hfill$\Box$

\subsection{A discrete ergodic decomposition for contractive Markov systems}\label{deds}

\begin{lemma}\label{evsl}
    Suppose $\M$ is positive.  Let $\mu\in P(\M)$ and $i\in N$ such that $\mu(K_i)>0$. If $j\in N$ is accessible from $i$, then $\mu(K_j)>0$.
\end{lemma}
{\it Proof.} 
Let $(e'_1,...,e'_n)$ be a path of $\M$ such that $i(e'_1)=i$ and $t(e'_n)=j$. Then
\begin{eqnarray*}
 \mu(K_j)&=&\int U^n1_{K_j}d\mu\\
&=&\sum\limits_{e_1,...,e_n}\int p_{e_1}p_{e_2}\circ w_{e_1}...p_{e_{n}}\circ w_{e_{n-1}}\circ...\circ w_{e_{1}}1_{K_j}\circ  w_{e_{n}}\circ...\circ w_{e_{1}}d\mu\\
&\geq&\int\limits_{K_i} p_{e'_1}p_{e'_2}\circ w_{e'_1}...p_{e'_{n}}\circ w_{e'_{n-1}}\circ...\circ w_{e'_{1}}1_{K_j}\circ  w_{e'_{n}}\circ...\circ w_{e'_{1}}d\mu\\
&>&0.
\end{eqnarray*}
\hfill$\Box$

\begin{Definition}
  For $f:K\longrightarrow\R$ and  $j\in N$, set
\[\Delta_j f(t):=\sup\limits_{x,y\in K_{j},d(x,y)\leq t} \left|f(x)-f(y)\right|\]
and, for $i\in N$ and $k\in\N$,
\[\Delta^{(k)}_i f(t):=\sup\left\{\Delta_j f(t)|\ j\mbox{ is accessible from }i\mbox{ by a path of length }k\right\}.\]
Let $0<\alpha<1$. We call $f$ $(\M,\alpha)$-{\it uniformly continuous} if and only if  for every $i\in N$ and $\beta>0$,
\[\lim\limits_{n\to\infty}\frac{1}{n}\sum\limits_{k=1}^n\Delta^{(k)}_i f\left(\beta a^{\alpha k}\right) = 0.\]
Note that for a bounded $f$, by Koopman-von Neumann Lemma (e.g. \cite{Pet}), this is equivalent to the existence of a set $B\subset\N$ of density zero (i.e.\\ $\lim_{n\to\infty}1/n\sum_{k=1}^n1_B(k)=0$) such that 
\[\lim\limits_{n\to\infty, n\notin B}\Delta^{(n)}_i f\left(\beta a^{\alpha n}\right) = 0.\]
Hence, the $(\M,\alpha)$-uniformly continuity of a bounded $f$ is equivalent to the existence of  $\gamma>0$ such that for every $i\in N$ and $\beta>0$,
\[\lim\limits_{n\to\infty}\frac{1}{n}\sum\limits_{k=1}^n\left(\Delta^{(k)}_i f(\beta a^{\alpha k})\right)^\gamma = 0.\]
\end{Definition}

\begin{lemma}\label{ashl}
    Suppose $\M$  is uniformly continuous, contractive with a contraction rate $0<a<1$ and satisfies Condition \ref{acC}. For each $e\in E$, let  $f_e:K\longrightarrow[0,\infty)$  be Borel-measurable such that  $\sum_{e\in E}p_ef_e$ is $(\M, \alpha)$-uniformly continuous for some $0<\alpha<1$ and $\sum_{e\in E}p_ef^2_e$ is bounded.  Let $x,y\in K_i$ for some $i\in N$. Then
\[\lim\limits_{n\to\infty}\frac{1}{n}\sum\limits_{k=1}^{n}\left(f_{\sigma_{k+1}}\circ w_{\sigma_k}\circ...\circ w_{\sigma_1}(x)-f_{\sigma_{k+1}}\circ w_{\sigma_k}\circ...\circ w_{\sigma_1}(y)\right)=0\ \ \ P_x\mbox{-a.e..}\]
\end{lemma}
{\it Proof.} 
We will use some ideas from \cite{Bre}.  For  $\sigma\in\S$ and $n\in\N$, let us abbreviate
 \[X^x_n(\sigma):=w_{\sigma_n}\circ w_{\sigma_{n-1}}\circ...\circ w_{\sigma_1}(x)\mbox{ and }X^x_0:=x,\]
 \[Y^x_{n}(\sigma):=f_{\sigma_{n}}\circ X^x_{n-1}(\sigma) \]
and
 \[Z^x_{n}:=Y^x_{n}-E_{P_x}\left(Y^x_{n}|\B_{n-1}\right)\mbox{ if }n\geq 2\mbox{ and }Z^x_{1} := 0.\]
We will apply the following result from \cite{Lo} p. 387. Let $Z_1,Z_2,...$ be a real-valued random process such that $E(Z_{n+1}|Z_{n},...,Z_{1})=0$ and $\sup_{n\in\N}E(Z^2_n)<\infty$ for all $n\in\N$. Then
\begin{eqnarray}\label{slln}
  \lim\limits_{n\to\infty}\frac{1}{n}\sum\limits_{k=1}^{n}Z_k=0\ \ \ \mbox{a.s.}
\end{eqnarray}
Note that 
\[E_{P_x}\left(Z^x_{n+1}|Z^x_{n},...,Z^x_{1}\right)=E_{P_x}\left(E_{P_x}\left(Z^x_{n+1}|\B_{n}\right)|Z^x_{n},...,Z^x_{1}\right)=0\]
$P_x$-a.e. for all $n\in\N$. Furthermore, one easily checks that
\begin{eqnarray}\label{evce}
  E_{P_x}\left(Y^x_{n+1}|\B_{n}\right)=\left(\sum\limits_{e\in E}p_ef_e\right)\circ X^x_n\ \ \ \ P_x\mbox{-a.e. for all }n\in\N.
\end{eqnarray}
Note that $\sum_{e\in E}p_ef_e\leq(\sum_{e\in E}p_ef_e^2)^{1/2}$. Therefore, by the  the hypothesis, there exists $0\leq\xi<\infty$ such that
\begin{eqnarray*}
&&\int \left(Z^x_{n}\right)^2dP_x\\
&=&\int \left(f_{\sigma_{n}}\circ  w_{\sigma_{n-1}}\circ...\circ w_{\sigma_1}(x)-\left(\sum\limits_{e\in E}p_ef_e\right)\circ  w_{\sigma_{n-1}}\circ...\circ w_{\sigma_1}(x)\right)^2dP_x(\sigma)\\
&=&\sum\limits_{e_1,...,e_n\in E}P_x\left(\L_1[e_1,...,e_{n-1}]\right)\left(p_{e_n}\left(f_{e_n}-\sum\limits_{e\in E}p_ef_e\right)^2\right)\circ w_{\sigma_{n-1}}\circ...\circ w_{\sigma_1}(x)\\
&=&U^{n-1}\left(\sum\limits_{e_n\in E}p_{e_n}\left(f_{e_n}-\sum\limits_{e\in E}p_ef_e\right)^2\right)(x)\\
&=&U^{n-1}\left(\sum\limits_{e\in E}p_ef_e^2\right)(x)-U^{n-1}\left(\left(\sum\limits_{e\in E}p_ef_e\right)^2\right)(x)\\
&\leq&\xi
\end{eqnarray*}
 for all $n\in\N$. Also, the same way,
\[\sup\limits_{n\in\N}\int \left(Z^y_{n}\right)^2dP_y<\infty.\]
Therefore, by \eqref{slln}, \eqref{evce} and the hypothesis,
\begin{eqnarray}\label{eecb}
&&\limsup\limits_{n\to\infty}\left|\frac{1}{n}\sum\limits_{k=1}^{n}\left(Y^x_k-Y^y_k\right)\right|\nonumber\\
&\leq&\limsup\limits_{n\to\infty}\frac{1}{n}\sum\limits_{k=1}^{n}\left|\left(\sum\limits_{e\in E}p_ef_e\right)\circ X^x_{k-1}-\left(\sum\limits_{e\in E}p_ef_e\right)\circ X^y_{k-1}\right|
\end{eqnarray}
$P_x$-a.e.. Furthermore, by the contraction on average property,
\[\int d(w_{\sigma_k}\circ...\circ w_{\sigma_1}x,w_{\sigma_k}\circ...\circ w_{\sigma_1}y)\ dP_x\leq a^kd(x,y)\]
for all $k\in\N$. Set $A^\alpha_{xyk}:=\{\sigma\in\Sp|\ d(w_{\sigma_k}\circ...\circ w_{\sigma_1}x,w_{\sigma_k}\circ...\circ w_{\sigma_1}y)>a^{\alpha k}d(x,y)\}$ for all $k\in\N$, and
\[A^\alpha_{xy}:=\bigcap\limits_{n=1}^\infty\bigcup\limits_{k=n}^\infty A^\alpha_{xyk}.\]
Then,  $P_x(A^\alpha_{xyk})\leq a^{(1-\alpha)k}$ for all $k\in\N$, and therefore, by the Borel-Cantelli argument, 
\begin{equation}\label{bszm}
   P_x(A^\alpha_{xy})=0.
\end{equation}
Note that for   $\sigma\in\Sp\setminus A^\alpha_{xy}$  there exists $m\in\N$ such that $d(w_{\sigma_k}\circ...\circ w_{\sigma_1}x,w_{\sigma_k}\circ...\circ w_{\sigma_1}y)\leq a^{\alpha k}d(x,y)$ for all $k\geq m$. Let $h:=\sum_{e\in E}p_ef_e$. Then, by the $(\M,\alpha)$-uniform equicontinuity and  boundedness of $h$,
\begin{eqnarray*}
&&\limsup\limits_{n\to\infty}\frac{1}{n}\sum\limits_{k=1}^{n}\left|h\circ X^x_{k-1}(\sigma)-h\circ X^y_{k-1}(\sigma)\right|\\
&\leq&\limsup\limits_{n\to\infty}\frac{1}{n}\sum\limits_{k=m+1}^{n}\Delta^{(k-1)}_i h\left(a^{\alpha (k-1)}d(x,y)\right)\\
&=&0.
\end{eqnarray*}
Thus, the assertion follows  by \eqref{bszm} and \eqref{eecb}.
\hfill$\Box$

\begin{Remark}\label{ekcr}
 An important special case for Lemma \ref{ashl} is when, for each $e\in E$, $f_e(x):=-\log p_e(x)$ if $x\in K_{i(e)}$ and $f_e(x):=0$ otherwise. In this case, the boundedness condition on $\sum_{e\in E}p_ef^2_e$ is always satisfied  if $E$ is finite. For, observe that  function
\[g(x):=-\sqrt{x}\log x\mbox{ if } x>0\mbox{ and }g(0):=0\]
is concave. Suppose $|E|=n$. Then
\[\frac{1}{n}\sum\limits_{e\in E}g(p_e)\leq g\left(\frac{1}{n}\right)=\frac{1}{\sqrt{n}}\log n,\]
and therefore,
\[\sum\limits_{e\in E}p_e\left(\log p_e\right)^2=\sum\limits_{e\in E}g^2\left(p_e\right)\leq\left(\sum\limits_{e\in E}g\left(p_e\right)\right)^2\leq n(\log n)^2.\]
However,  as Lemma \ref{edl} shows, ACC is  actually the only hypothesis which is needed in this case.
\end{Remark}

\begin{lemma}\label{edl}
    Suppose $\M$  satisfies Condition \ref{acC}. For each $e\in E$, let $f_e(z):=-\log p_e(z)$ if $z\in K_{i(e)}$ and $f_e(z):=0$ otherwise.
  Let $x,y\in K_i$ for some $i\in N$. Then
\[\lim\limits_{n\to\infty}\frac{1}{n}\sum\limits_{k=1}^{n}\left(f_{\sigma_{k+1}}\circ w_{\sigma_k}\circ...\circ w_{\sigma_1}(x)-f_{\sigma_{k+1}}\circ w_{\sigma_k}\circ...\circ w_{\sigma_1}(y)\right)=0\ \ \ P_x\mbox{-a.e..}\]
\end{lemma}
{\it Proof.} 
Observe that,  since $P_x(\Sp_G)=1$, by the hypothesis,
\begin{eqnarray*}
 &&\frac{1}{n}\sum\limits_{k=1}^{n}\left(f_{\sigma_{k+1}}\circ w_{\sigma_k}\circ...\circ w_{\sigma_1}(x)-f_{\sigma_{k+1}}\circ w_{\sigma_k}\circ...\circ w_{\sigma_1}(y)\right)\\
&=&\frac{1}{n}\log\frac{P_y\left(\L_1[\sigma_1,...,\sigma_{n+1}]^+\right)}{P_x\left(\L_1[\sigma_1,...,\sigma_{n+1}]^+\right)}+\frac{1}{n}\log\frac{p_{\sigma_1}(x)}{p_{\sigma_1}(y)}\\
&=&\frac{1}{n}\log Z^{yx}_{n+1}(\sigma)+\frac{1}{n}\log\frac{p_{\sigma_1}(x)}{p_{\sigma_1}(y)},
\end{eqnarray*}
and $ Z^{yx}_{n}(\sigma)=1/Z^{xy}_{n}(\sigma)$ for $P_x$-a.e. and $P_y$-a.e. $\sigma\in\Sp$ and all $n\in\N$. Hence, 
$ Z^{yx}_\infty=1/Z^{xy}_\infty$ $P_x$-a.e. and $P_y$-a.e., with the definitions $ 1/0:=\infty$ and $1/\infty=0$. By Theorem \ref{acet}, $P_y\{ Z^{yx}_\infty<\infty\}=1$ and $P_x\{Z^{xy}_\infty<\infty\}=1$. Hence, 
\[P_x\left\{0< Z^{yx}_\infty<\infty\right\}=1.\]
That is
\[P_x\left\{\left|\log Z^{yx}_\infty\right|<\infty\right\}=1.\]
Therefore,
\[P_x\left\{\lim\limits_{n\to\infty}\frac{1}{n}\left|\log Z^{yx}_{n+1}\right|=0\right\}=1,\]
as desired. \hfill$\Box$

\begin{theo}\label{uim}
  Suppose $\M$  is uniformly continuous, contractive with a contraction rate $0<a<1$ and satisfies Condition \ref{acC}. Let $M\in E(\M)$ and $I_M:=\{i\in I|\ F(M)(C_i)>0\}$.   For each $e\in E$, let $f_e:K\longrightarrow[0,\infty)$  be Borel-measurable such that\\
(a) $\sum_{e\in E}p_ef_e$ is $(\M, \alpha)$-uniformly continuous for some $0<\alpha<1$ and $\sum_{e\in E}p_ef^2_e$ is bounded, or\\
(b)   $f_e(z):=-\log p_e(z)$ if $z\in K_{i(e)}$ and $f_e(z):=0$ otherwise, and $-\sum_{e\in E}p_e\log p_e$ is bounded.\\
 Then the following holds true.\\
(i) $I_M$ is not empty. For every $i\in I_M$,  there exists a unique $\Lambda_i\in E(\M)$ with $F(\Lambda_i)(C_i)=1$.  For every $i\in I_M$, $\Lambda_i$ is ergodic, $F^{-1}(C_i)\in\I_M$ and $\{\mu\in P(\M)|\ \mu(C_i)=1\} = \{F(\Lambda_i)\}$. For every $B\in\BS$ and $Q\in\I_M$,
\[M\left(B\cap Q\right)=\sum\limits_{i\in I_M}M\left(F^{-1}(C_i)\cap Q\right)\Lambda_i\left(B\right).\]
(ii)  For every $x\in \bigcup_{i\in I_M}C_i$,
\[\frac{1}{n}\sum\limits_{k=1}^{n}f_{\sigma_{k+1}}\circ{w_{\sigma_k}\circ...\circ w_{\sigma_1}(x)}\to
\sum\limits_{i\in I_M}1_{C_i}(x)\sum\limits_{e\in E}\int\limits_{K_{i(e)}} p_ef_edF(\Lambda_i)\]
 for $P_x$-a.e. $\sigma\in\Sigma^+$. \\
(iii) $I_M\subset I_G$. If $\M$ has a single generating class, then $M$ is ergodic and $E(\M)=\{M\}$.
\end{theo}
{\it Proof.}  
(i) Let $M=\int\Lambda_\omega dM(\omega)$ be the ergodic decomposition of $M$ given by Theorem \ref{ted}. Let $f_\infty := 0$ and $v(\sigma):=f_{\sigma_1}\circ F(\sigma)$ for all $\sigma\in\S$. Then, by the hypothesis,
\begin{eqnarray*}
 \infty&>&\sup_{x\in K}\left(\sum\limits_{e\in E}p_e(x)f^2_e(x)\right)^{\frac{1}{2}}\geq\sup_{x\in K}\sum\limits_{e\in E}p_e(x)f_e(x)\geq\int\sum\limits_{e\in E}p_e\circ Ff_e\circ FdM\\
&=&\sum\limits_{e\in E}\int E_M\left(1_{_1[e]}|\F\right)f_e\circ FdM=\sum\limits_{e\in E}\int 1_{_1[e]}f_e\circ FdM=\int |v|dM.
\end{eqnarray*}
 Hence, $v\in \mathcal{L}^1(M)$ in both cases, (a) and (b). Therefore, by Corollary  \ref{esdc} (i), there exists $X\in\BS$ with $M(X)=1$ such that  $\Lambda_\omega\in E(\M)$ and $v\in\mathcal{L}^1(\Lambda_\omega)$  and for all $\omega\in X$. 

Now, let $\omega\in X$. Note that $F\circ S^k(\sigma) = \bar w_{\sigma_k}\circ...\circ \bar w_{\sigma_1} (F(\sigma))$ for all $k\in\N$ and $\sigma\in D$. Then, by Birkhoff's ergodic theorem,
\[\frac{1}{n}\sum\limits_{k=1}^{n}f_{\sigma_{k+1}}\circ\bar w_{\sigma_k}\circ...\circ\bar w_{\sigma_1}(F(\sigma))=\frac{1}{n}\sum\limits_{k=1}^{n}v\circ S^k(\sigma)\to\int v\ d\Lambda_\omega\ \ \ \mbox{ for } \Lambda_\omega\mbox{-a.e. }\sigma\in\S.\]
Set 
\[\bar f_n(x,\sigma):=\frac{1}{n}\sum\limits_{k=1}^{n}f_{\sigma_{k+1}}\circ\bar w_{\sigma_k}\circ...\circ\bar w_{\sigma_1}(x)\ \ \ \mbox{ for all }(x,\sigma)\in K\times\Sp\mbox{ and }n\in\N,\]
\begin{eqnarray*}
  \eta:\Sigma&\longrightarrow& K\times\Sigma^+\\
       \sigma&\longmapsto&(F(\sigma),(\sigma_1,\sigma_2,...)),
\end{eqnarray*}
and
\[Q_\omega:=\left\{\sigma\in\S|\ \bar f_n\circ\eta(\sigma)\to\int v\ d\Lambda_\omega\right\}.\]
Then $\Lambda_\omega(Q_\omega)=1$, and
\[\bar f_n(x,\sigma)\to\int v\ d\Lambda_\omega\ \ \ \mbox{ for all }(x,\sigma)\in\eta(Q_\omega).\]
Since $\Lambda_\omega\in E(\M)$, by Lemma 4 (ii) in \cite{Wer11}, $\eta(\Lambda_\omega)=\tilde\phi(F(\Lambda_\omega))$. Hence, by the definition of $\tilde\phi$,
\begin{eqnarray*}
 1&=&\Lambda_\omega(Q_\omega)\\
&\leq&\eta(\Lambda_\omega)(\eta(Q_\omega))\\
&=&\tilde\phi(F(\Lambda_\omega))(\eta(Q_\omega))\\
&=&\int P_x\left(\left\{\sigma\in\Sp|\ (x,\sigma)\in\eta(Q_\omega)\right\}\right)dF(\Lambda_\omega)(x)\\
 &=&\sum\limits_{i\in S_\omega} \int\limits_{K_i}
P_x\left(\left\{\sigma\in\Sp|\ (x,\sigma)\in\eta(Q_\omega)\right\}\right)dF(\Lambda_\omega)(x)
\end{eqnarray*}
where $S_\omega:=\{i\in N|\ F(\Lambda_\omega)(K_i)>0\}$. Therefore, for every $i\in S_\omega$ there exists $y_i\in K_i$ such that 
\[P_{y_i}\left(\left\{\sigma\in\Sp|\ (y_i,\sigma)\in\eta(Q_\omega)\right\}\right)=1.\]
Let $j\in S_\omega$ and set $\Omega_j:=\{\sigma\in\Sp|\ (y_j,\sigma)\in\eta(Q_\omega)\}$. Then,   for every $\sigma\in\Omega_j$,
\[\bar f_n(y_j,\sigma)\to\int v\ d\Lambda_\omega.\]
Set
\[ f_n(x,\sigma):=\frac{1}{n}\sum\limits_{k=1}^{n}f_{\sigma_{k+1}}\circ w_{\sigma_k}\circ...\circ w_{\sigma_1}(x)\ \ \ \mbox{ for all }(x,\sigma)\in K\times\Sp\mbox{ and }n\in\N.\]
As $P_{y_j}(\{\sigma\in\Sp|\ y_j\notin K_{i(\sigma_1)}$ or $\exists k\in\mathbb{N}$ s.t. $i(\sigma_{k+1})\neq t(\sigma_{k})\}) = 0$, without a loss of generality, we can assume that for every $\sigma\in\Omega_j$,
\[ f_n(y_j,\sigma)\to\int v\ d\Lambda_\omega\]
and $P_{y_j}(\Omega_j) = 1.$
Let $y\in K_j$. Since $P_y\ll P_{y_j}$, $P_y(\Omega_j) = 1$,  by Lemma \ref{ashl} in case (a) and by Lemma \ref{edl} in case (b),
\begin{eqnarray*}
 \limsup\limits_{n\to\infty}\left|f_n(y,\sigma)-\int v\ d\Lambda_\omega\right|\leq \limsup\limits_{n\to\infty}\left|f_n(y,\sigma)-f_n(y_j,\sigma)\right|=0
\end{eqnarray*}
for $P_y$-a.e. $\sigma\in\Sp$. Thus, for every  $y\in\bigcup_{j\in S_{\omega}}K_j$,
\begin{equation}\label{eec}
f_n(y,\sigma)\to\int\sum_{e\in E}1_{_1[e]}f_e\circ Fd\Lambda_\omega=\sum_{e\in E}\int p_ef_edF(\Lambda_\omega)\ \ \ \mbox{ for }P_y\mbox{-a.e }\sigma\in\Sp.
\end{equation}
In particular, for every $y\in\bigcup_{j\in S_{\omega}}K_j$ and $i\in N$,
\begin{equation}\label{ecc}
   \frac{1}{n}\sum\limits_{k=1}^{n}1_{K_i}\circ w_{\sigma_k}\circ...\circ w_{\sigma_1}(y)\to F(\Lambda_\omega)(K_i)\ \ \ \mbox{ for }P_y\mbox{-a.e }\sigma\in\Sp.
\end{equation}
This implies that every $i\in S_\omega$ is essential. Fix $i_0\in S_\omega$, and let $c_{i_0}\subset N$ denote the communication class containing $i_0$. Let $j\in  S_\omega$. Then, by  \eqref{ecc}, $j$ is accessible from $i_0$ and vice versa. Hence, $j\in  c_{i_0}$. Thus $S_\omega\subset c_{i_0}$. Let $i\in c_{i_0}$. Then there exists $j\in  S_\omega$ which communicates with $i$. By Proposition 1 in \cite{Wer11}, $F(\Lambda_\omega)\in P(\M)$, and therefore $F(\Lambda_\omega)(K_i)>0$ by Lemma \ref{evsl}. Hence, $i\in S_\omega$. That is $c_{i_0}\subset  S_\omega$. Thus $S_\omega = c_{i_0}$. This defines a map
\begin{eqnarray*}
       \theta :X&\longrightarrow& I\\
       \omega&\longmapsto&\mbox{index of the communication class } S_\omega.
\end{eqnarray*}
As for every $j\in I$, $\theta^{-1}\{j\}=\{\omega\in X|\ S_\omega=c_j\}=\{\omega\in X|\ \Lambda_\omega(F^{-1}(C_j))=1\}$, $\theta$ is $\I_M$-Borel-measurable, and
\begin{eqnarray}\label{pame}
      M\left(F^{-1}(C_j)\right)&=&\int\limits_{X}\Lambda_\omega\left(F^{-1}(C_j)\right)dM(\omega)\nonumber\\
 &=&\sum\limits_{i\in I}\int\limits_{\theta^{-1}\{i\}}\Lambda_\omega\left(F^{-1}(C_j)\right)dM(\omega)\nonumber\\
&=&\sum\limits_{i\in I}\int\limits_{\left\{\omega\in X|\ \Lambda_\omega(F^{-1}(C_i))=1\right\}}\Lambda_\omega\left(F^{-1}(C_j)\right)dM(\omega)\nonumber\\
&=&M\left(\{\omega\in X|\ \Lambda_\omega(F^{-1}(C_j))=1\}\right)\nonumber\\
&=&M\left(\theta^{-1}\{j\}\right).
\end{eqnarray}
Hence, $ I_M$ is not empty. Let $i\in I_M$. Then  $\theta^{-1}\{i\}$ is not empty.    Fix $\omega^i\in \theta^{-1}\{i\}$.  Then, by \eqref{eec}, for every $y\in C_i$,
\begin{equation}\label{edec}
f_n(y,\sigma)\to\sum_{e\in E}\int p_ef_edF(\Lambda_{\omega^i})\ \ \ \mbox{ for }P_y\mbox{-a.e }\sigma\in\Sp.
\end{equation}
Let $f:K\longrightarrow\R$ be uniformly continuous and bounded. Then,  for every $y\in C_i$,
\begin{equation}\label{fec}
   \frac{1}{n}\sum\limits_{k=1}^{n}f\circ w_{\sigma_k}\circ...\circ w_{\sigma_1}(y)\to\int f dF(\Lambda_{\omega^i})\ \ \ \mbox{ for }P_y\mbox{-a.e }\sigma\in\Sp.
\end{equation}
Hence, the integration of \eqref{fec} by $P_y$ implies by Lebesgue's Dominated Convergence Theorem that,  for every $y\in C_i$,
\begin{equation}\label{eoc}
 \frac{1}{n}\sum\limits_{k=1}^{n}U^kf(y)\to\int f dF(\Lambda_{\omega^i}).
\end{equation}
Clearly, $F(\Lambda_{\omega^i})\in P(\M)$ with $F(\Lambda_{\omega^i})(C_i)=1$. Suppose there is another $\mu_i\in P(\M)$  with $\mu_i(C_i)=1$. Then the integration of \eqref{eoc} by $\mu_i$  implies by Lebesgue's Dominated Convergence Theorem that
\begin{equation}\label{edim}
      \int fd\mu_i=\int  f dF(\Lambda_{\omega^i}).
\end{equation}
That is 
\begin{equation}\label{ueim}
      \mu_i=F(\Lambda_{\omega^i}),
\end{equation}
 since the set of all bounded uniformly continuous functions separates the measures. Hence, by Proposition 1 and Lemma 4 (i) in \cite{Wer11}, for every $\Lambda\in E(\M)$ with $F(\Lambda)(C_i)=1$,
\begin{equation}\label{uees}
      \Lambda=\Phi(F(\Lambda))=\Phi(F(\Lambda_{\omega^i})) = \Lambda_{\omega^i}.
\end{equation}
In particular, it follows that $\Lambda_{\omega}=\Lambda_{\omega^i}$ for all $\omega\in \theta^{-1}\{i\}$ and $i\in I_M$. Therefore, for every $Q\in\I_M$ and $B\in\BS$,
\[M\left(B\cap Q\right)=\int\limits_Q\Lambda_\omega\left(B\right)dM(\omega)=\sum\limits_{i\in I_M}M\left(\theta^{-1}\{i\}\cap Q\right)\Lambda_{\omega^i}\left(B\right).\]
In particular, for every $i\in I_M$,
\begin{equation}\label{sicc}
      M\left(F^{-1}(C_i)\cap \theta^{-1}\{i\}\right)=M\left(\theta^{-1}\{i\}\right).
\end{equation}
Together with \eqref{pame}, this implies that $M\left(F^{-1}(C_i)\Delta \theta^{-1}\{i\}\right)=0$ for all $i\in I_M$. Thus, for every $i\in I_M$, $F^{-1}(C_i)\in\I_M$, since $\theta^{-1}\{i\}\in\I_M$, and for every 
 $Q\in\I_M$ and $B\in\BS$,
\begin{equation}\label{fede}
     M\left(B\cap Q\right)=\sum\limits_{i\in I_M}M\left(F^{-1}(C_i)\cap Q\right)\Lambda_{\omega^i}\left(B\right).
\end{equation}
 This competes the proof of  (i), with $\Lambda_{\omega^i}$ for $\Lambda_{i}$.

(ii)  follows immediately  from \eqref{edec} with $\Lambda_{\omega^i}$ for $\Lambda_{i}$.  

(iii)  By  \eqref{edec}, for every $j\in I_M$, $y\in C_j$ and $i\in N$,
\begin{equation}
   \frac{1}{n}\sum\limits_{k=1}^{n}1_{K_i}\circ w_{\sigma_k}\circ...\circ w_{\sigma_1}(y)\to F(\Lambda_{\omega_j})(K_i)\ \ \ \mbox{ for }P_y\mbox{-a.e }\sigma\in\Sp.
\end{equation}
Hence, the integration  by $P_y$ implies that each communication class in $I_M$ is generating. Let $C_g$ be the singe generating class of $\M$. Then, by \eqref{fede}, $M$ is ergodic with $M(F^{-1}(C_g))=1$. Let $M'\in E(M)$. Then, since $M$ was an arbitrary member of $E(\M)$, $M'(F^{-1}(C_g))=1$, and therefore, by  \eqref{uees}, $M=M'$. 
\hfill$\Box$

\subsection{Controllable invariant measures}\label{cim}

Now, we are going to apply Theorem \ref{uim} for a description of all invariant Borel probability measures of such $\M$ in the non-degenerate case by means of the generating communication classes, through the usage of the results on the existence of the equilibrium states and their bijective correspondence to the invariant Borel probability measures  which were obtained in \cite{Wer11}.

\begin{Definition}
  We call $\mu\in P(\M)$ ergodic if and only if for every $\lambda\in P(\M)$, $\lambda\ll\mu$ implies $\lambda=\mu$. 
\end{Definition}
Clearly, a unique invariant Borel probability measure is ergodic.

\begin{Condition}\label{dtc}[Conditions for non-emptiness of $\tilde E(\M)$ for a contractive uniformly continuous $\M$ with a generating communication class, by Theorem 5 in \cite{Wer11}]
\newline (i) 
\begin{equation}\label{cgc}
       \sup\limits_{i\in N}\sup\limits_{x\in K_i}\sum\limits_{e\in E,\ i(e) = i}p_e(x)d\left( w_e(x_{i(e)}), x_{t(e)}\right)<\infty,
\end{equation}  
(ii) For every $i\in N$,
\begin{equation}\label{dmcc}
       \sum\limits_{e\in E, i(e)=i}\sup\limits_{x\in K_i}p_e(x)<\infty.
 \end{equation} 
 $\M$ is said to have a {\it dominating Markov chain} if (ii) is satisfied.
\end{Condition}

\begin{theo}\label{esgc}
 Suppose $\M$  is uniformly continuous, non-degenerate, contractive and satisfies Conditions \ref{acC} and \ref{dtc}.   Then the following holds true.\\
(i) For every $i\in I_G$, there exists a unique $\Lambda_i\in E(\M)$ such that $F(\Lambda_i)(C_i)=1$. For every $i\in I_G$,  $\Lambda_i$ is ergodic.\\
(ii) For every $i\in I_G$, there exists a unique $\lambda_i\in P(\M)$ such that $\lambda_i(C_i)=1$, and therefore, $\lambda_i$ is ergodic.   For each $e\in E$, let $f_e:K\longrightarrow[0,\infty)$  be Borel-measurable  such that \\
(a) $\sum_{e\in E}p_ef_e$ is $(\M, \alpha)$-uniformly continuous for some $0<\alpha<1$ and $\sum_{e\in E}p_ef^2_e$ is bounded, or\\
(b)   $f_e(z):=-\log p_e(z)$ if $z\in K_{i(e)}$ and $f_e(z):=0$ otherwise, and $-\sum_{e\in E}p_e\log p_e$ is bounded.\\
  Then, for every $i\in I_G$ and $x\in C_i$,
\[\frac{1}{n}\sum\limits_{k=1}^{n}f_{\sigma_{k+1}}\circ{w_{\sigma_k}\circ...\circ w_{\sigma_1}(x)}\to
\sum\limits_{e\in E}\int\limits_{K_{i(e)}} p_ef_ed\lambda_i\ \ \ \mbox{ for }P_x\mbox{-a.e. }\sigma\in\Sigma^+.\]
(iii) For every $\mu\in P(\M)$,
\[\mu\left(B\right)=\sum\limits_{i\in I_G}\mu\left(C_i\right)\lambda_i\left(B\right)\ \ \ \mbox{ for all }B\in\B(K).\]
(iv) $\M$ has a unique invariant Borel probability measure if and only if it has a single generating communication class.
\end{theo}
{\it Proof.} 
 (i) Let $i'\in I_G$ and $x_0\in K_i\subset C_{i'}$ be a generating point. Then,  by Theorem 5 (i) in \cite{Wer11}, there exists $\Lambda_{x_0}\in\tilde E(\M)$ such that for every $j\in N$  with $\Lambda_{x_0}(T_j)>0$ there exists a path from $i$  to  $j$. Hence $\{j\in N|\ \Lambda_{x_0}(T_j)>0\}\subset c_{i'}$. Since $\M$ is non-degenerate, by Theorem 2 in \cite{Wer11}, $\Lambda_{x_0}\in E(\M)$, and therefore, for every $j\in N$,
\[\Lambda_{x_0}(T_j) = \sum\limits_{e\in E,i(e)=j}\int1_{_1[e]}d\Lambda_{x_0}= \sum\limits_{e\in E,i(e)=j}\int p_e\circ Fd\Lambda_{x_0}=\Lambda_{x_0}\left(F^{-1}(K_j)\right).\]
Hence, $\Lambda_{x_0}\left(F^{-1}(C_{i'})\right)=1$. By Theorem \ref{uim} (i), $\Lambda_{x_0}$ is unique in $E(\M)$ with such a property and is ergodic.

(ii) Let $i\in I_G$. Set $\lambda_i:=F(\Lambda_i)$, where $\Lambda_i$ is given by (i). Then $\lambda_i\in P(\M)$ with $\lambda_i(C_i)=1$. By Theorem \ref{uim} (i), it is unique with such a property, and therefore, ergodic. The rest of (ii) follows by Theorem \ref{uim} (ii).

(iii) Let $\mu\in P(\M)$. Then, by the non-degeneracy of $\M$, $\Phi(\mu)\in E(\M)$, by Theorem 5 (ii) in \cite{Wer11} and Theorem 2 in \cite{Wer11}. Hence, by Theorem \ref{uim} (i) and  (iii), for every $B\in \BS$,
\[\Phi(\mu)(B)=\sum\limits_{i\in I_{\Phi(\mu)}}\Phi(\mu)\left(F^{-1}(C_i)\right)\Lambda_i\left(B\right)=\sum\limits_{i\in I_G}\Phi(\mu)\left(F^{-1}(C_i)\right)\Lambda_i\left(B\right).\]
By  Corollary 1 (ii) in \cite{Wer11}, $\Phi$ is the inverse of $F:E(\M)\longrightarrow P(\M)$. Therefore, for every $B\in\B(K)$,
\[\mu(B)=\Phi(\mu)(F^{-1}(B))=\sum\limits_{i\in I_G}\mu(C_i)\Lambda_i\left(F^{-1}(B)\right)=\sum\limits_{i\in I_G}\mu(C_i)\lambda_i(B).\]

(iv) The 'if' part follows by (i), Theorem \ref{uim} (iii) and Corollary 1 (ii) in \cite{Wer11}. Now, suppose $P(\M)=\{\mu\}$. Then, as above, $\Phi(\mu)\in E(\M)$, and by Theorem  \ref{uim} (i) and (iii), $I_G$ is not empty. Let $i\in I_G$. By (ii), there exists unique $\lambda_i\in P(\M)$ such that $\lambda_i(C_{i})=1$. Thus $\mu(C_{i})=1$. Therefore, there can be only one generating communication class.
 \hfill$\Box$

The following is a straightforward application of Theorem \ref{esgc} which allows to compute the entropy of processes generated by $\M$.
\begin{cor}\label{Ec}
 Suppose $\M$  is uniformly continuous, non-degenerate, contractive and satisfies Conditions \ref{acC} and \ref{dtc}. Suppose    $-\sum_{e\in E}p_e\log p_e$ is bounded.   Then  for every $i\in I_G$ there exists an ergodic $\lambda_i\in P(\M)$ such that for  every $x\in C_i$,
\[\lim\limits_{n\to\infty}\frac{1}{n}\log P_x( _1[\sigma_1,...,\sigma_n])=\sum\limits_{e\in E}\int\limits_{K_{i(e)}} p_e\log p_ed\lambda_i\]
for $P_x$-a.e.  $\sigma\in\Sigma^+$.
\end{cor}
{\it Proof.} 
By Theorem \ref{esgc} (ii)  case (b),  for every $i\in I_G$ and $x\in C_i$,
\begin{eqnarray*}
&&-\lim\limits_{n\to\infty}\frac{1}{n}\log [p_{\sigma_1}(x)p_{\sigma_2}\circ w_{\sigma_1}(x)...
p_{\sigma_n}\circ w_{\sigma_{n-1}}\circ...\circ
w_{\sigma_1}(x)]\\
&=&-\lim\limits_{n\to\infty}\frac{n-1}{n}\frac{1}{n-1}\sum\limits_{k=1}^{n-1}\log p_{\sigma_{k+1}}\circ{w_{\sigma_k}\circ...\circ w_{\sigma_1}(x)}-\lim\limits_{n\to\infty}\frac{1}{n}\log p_{\sigma_{1}}(x)\\
&=&-\sum\limits_{e\in E}\int\limits_{K_{i(e)}} p_e\log p_ed\lambda_i
\end{eqnarray*}
 for $P_x$-a.e.  $\sigma\in\Sigma^+$, as desired.  \hfill$\Box$

\subsection{Some sufficient conditions for the absolute continuity condition}\label{scac}

\subsubsection{Connecting refinements of a Markov system}

Often one can choose several Markov systems associated with a random dynamical system, in particular, such which form a tree with respect to the refinement of their Markov partitions. This can be very useful in establishing the validity of ACC for such Markov systems,  e.g. see Example \ref{se} for an illustration. 

The following definition of a refiniment of $\M$ was given in \cite{Wer11}.
\begin{Definition}
    We call a Markov system  $\M^r:=(K^r_{i(e)},w^r_e,p^r_e)_{e\in E^r}$ a {\it refinement} of $\M$ if and only if partition 
$\{K^r_{i(e)}\}_{e\in E^r}$ refines  partition $\{K_{i(e)}\}_{e\in E}$ (i.e. each $K_i$ is a union of some $K^r_j$'s) and there is a surjective map $r:E^r\longrightarrow E$ such that
$w_{r(e)}|_{K^r_{i(e)}} = w^r_{e}|_{K^r_{i(e)}}$ and $p_{r(e)}|_{K^r_{i(e)}} = p^r_{e}|_{K^r_{i(e)}}$ for all $e\in E^r$ (we use the same letters for maps $i,t:E^r\longrightarrow N^r$).  $r$  is called the {\it refinement map}. 
\end{Definition}

\begin{Definition}
We call  refinements $\M^1$ and $\M^2$ of $\M$ {\it connecting} if and only if for every $x,y\in K_i$ and $i\in N$ there exist $j_1,...,j_n\in N^1$ and $j'_1,...,j'_n\in N^2$ such that $x\in K^1_{j_1}\mbox{ and }y\in K^2_{j'_n}$, and
\begin{equation}\label{cpc}
   K^1_{j_i}\cap K^2_{j'_i}\ne\emptyset\mbox{ and }K^1_{j_i}\cap K^2_{j'_{i+1}}\ne\emptyset\mbox{ for all }i=1,..,n-1.
\end{equation}
\end{Definition}

\begin{lemma}\label{crl}
Let $\M^1$ and $\M^2$ be connecting refinements of $\M$. If $\M^1$ and $\M^2$ satisfy Condition \ref{acC}, then $\M$ satisfies Condition \ref{acC} also. 
\end{lemma}
{\it Proof.} Let $x,y\in K_j$ for some $j\in N$. Since $\M^1$ and $\M^2$ are connecting refinements of $\M$, there exist $j_1,...,j_n\in N^1$ and $j'_1,...,j'_n\in N^2$ such that $x\in K^1_{j_1}\mbox{ and }y\in K^2_{j'_n}$ and  \eqref{cpc} is satisfied. Then, applying Lemma 1 (iv) in \cite{Wer11} to refinement $\M^1$ implies that $P_x$ is equivalent to $P_z$ for all $z\in  K^1_{j_1}$. Since $K^1_{j_1}\cap K^2_{j'_{2}}\ne\emptyset$ there exists  $z_2\in K^2_{j'_{2}}$ such that $P_x$ is equivalent  to $P_{z_2}$. Therefore, applying Lemma 1 (iv) in \cite{Wer11} to refinement $\M^2$ implies that $P_x$ is equivalent to $P_z$ for all $z\in  K^2_{j'_{2}}$. The same way, since $ K^1_{j_2}\cap K^2_{j'_2}\ne\emptyset$, $P_x$ is equivalent to $P_z$ for all $z\in  K^1_{j_{2}}$. Thus,  repeating the argument $n-1$-times, implies that $P_x$ is equivalent to $P_z$ for all $z\in  K^2_{j'_{n}}$. Thus, in particular, $P_x$ is equivalent to $P_y$.
\hfill$\Box$

\subsubsection{Relative entropy}

\begin{Definition}
Suppose $\M$ is positive. For $i\in N$, $x,y\in K_i$ and $n\in\N$, set
 \[H_n(P_x|P_y):=\int Z^{xy}_n\log Z^{xy}_ndP_y\]
with the continuous extension $0\log 0:=0$. It is well known that $0\leq H_n(P_x|P_y)\leq H_{n+1}(P_x|P_y)$ for all $n\in\N$ (e.g. see  \cite{S}, p.78). 
Set $H(P_x|P_y):=\lim_{n\to\infty}H_n(P_x|P_y)$. It is called the {\it relative entropy}, {\it Kullback-Leibler entropy} or {\it Kullback-Leibler divergence} of measures.
\end{Definition}
A well-known sufficient condition for the absolute continuity is the following fact.
\begin{theo}\label{klfc}
    Suppose $\M$ is positive. Let $x,y\in K_i$ for some $i\in N$. If $H(P_x|P_y)<\infty$, then $P_x\ll P_y$.
\end{theo}
{\it Proof.} One can use a deeper result useful for checking the uniform integrability due to   Ch.-J. de la Vall\'{e}e Poussin, e.g.  Theorem 4.5.9 in \cite{Bog} Vol. 1 (see also Example 4.5.10 there) or \cite{Doob} and then conclude the assertion by Theorem \ref{acet}, or observe that, by the definition of $Z^{xy}_n$, for every $n\in\N$,
\[H_n(P_x|P_y)=\int\log  Z^{xy}_ndP_x.\]
Hence, the hypothesis implies that
\[\sup\limits_{n\in\N}P_x\left\{\log Z^{xy}_n>c\right\}\to 0\mbox{ as }c\to\infty.\]
Thus, the assertion follows by Theorem \ref{acet}.
\hfill$\Box$

\subsubsection{Square summability of variation of probability functions}

 In this subsection, we give a sufficient condition for the finiteness of the relative entropy in terms of the variation of the probability functions in the case of finitely many $e\in E$ with $i(e)=j$ for each $j\in N$.
 It is called the {\it square summability of variation}. It has been used by A. Johansson and A. \"{O}berg in \cite{JO}, where they proved the uniqueness of the $g$-measure for a $g$-function satisfying this condition and the boundedness away from zero.  N. Berger, Ch. Hoffman and V. Sidoravicius have shown in \cite{BHS} that the condition of the square summability of variation is tight, in the sense that for any $\epsilon>0$ there exists a $g$-function  with a summable variation to the power $2+\epsilon$ which has several $g$-measures. The reader is referred to \cite{JO} for a discussion on the relation between the Johansson-\"{O}berg condition and the Berbee condition \cite{Ber}, and to \cite{JOP} and \cite{JOP2}  for the latest on the weakening of the continuity of $g$-functions with unique $g$-measures (note that the $g$-functions associated with the random dynamical systems considered in this article are not continuous even in the case with an open Markov partition).
 
\begin{Definition}
  $f:K\longrightarrow\mathbb{R}$ is  said to have a {\it square summable variation} if and only if
for any $c >0$
 \[\int_0^c\frac{\phi^2(t)}{t}dt<\infty\]
  where $\phi$ is {\it the modulus of uniform continuity} of $f$, i.e.
     \[\phi(t):=\sup\{|f(x)-f(y)|:d(x,y)\leq t,\ x,y\in X\}\]
or equivalently, for any $b>0$ and $0<c<1$,
 \[\sum\limits_{n=0}^\infty\phi^2\left(bc^n\right)<\infty,\] which is obviously a weaker condition than
the {\it Dini-continuity (summability of variation)} and stronger than the uniform continuity.
\end{Definition}
 
\begin{Definition}
   Suppose $\M$ is positive and contractive with a contraction rate $0<a<1$. For $e\in E$, let $\phi_e$ be the modulus of uniform continuity of $p_e|_{K_{i(e)}}$.  Let $j\in N$ and $x,y\in K_i$ for some $i\in N$. For $0\leq \beta\leq 1$ and $k>0$, set  
\[A_{xy}^{(\beta k)}:=\left\{\sigma\in\Sigma^+_G:\ d(w_{\sigma_k}\circ...\circ w_{\sigma_1}x,w_{\sigma_k}\circ...\circ w_{\sigma_1}y)>a^{\beta k} d(y,x)\right\},\]
\[\ell_j:=\sup\limits_{z\in K_j}\sum\limits_{e\in E, i(e)=j}\frac{1}{p_e(z)},\]
\[\phi^{(\beta k)}_{xyj}:=\sup\limits_{z\in K_j}\sum\limits_{e\in E, i(e)=j}\frac{\phi^2_e\left(a^{\beta k}d(x,y)\right)}{p_e(z)}\]
 and
 \[B_{\beta xy}:=\sum_{k=0}^\infty\sum\limits_{j\in N}\left(\ell_j\int\limits_{A_{xy}^{(\beta k)}}1_{K_j}\circ w_{\sigma_k}\circ...\circ w_{\sigma_1}(x)dP_x(\sigma)+\phi^{(\beta k)}_{xyj}{U^*}^k\delta_x(K_j)\right)\]  
with $A_{xy}^{(\beta 0)}:=\Sigma^+_G$ if at least one of $\{w_e\}_{e\in E}$ is not a contraction and $A_{xy}^{(\beta 0)}:=\emptyset$ otherwise.
\end{Definition}

\begin{lemma}\label{acl}
 Suppose $\M$ is positive and contractive with a contraction rate $0<a<1$. Let $0\leq \beta\leq 1$ and $x,y\in K_i$ for some $i\in N$.
Then \[H(P_x|P_y)\leq\left(B_{\beta xy}\right)^\frac{1}{2}+B_{\beta xy}.\]
\end{lemma}
{\it Proof.} We adapt a part of the proof from \cite{JO}. Fix $j\in N$ and $x,y\in K_j$.
    Let us abbreviate
\[p_i^y(\sigma):=p_{\sigma_i}(w_{\sigma_{i-1}}\circ...\circ w_{\sigma_1}y)\]
for $i\geq 2$ and $p_1^y(\sigma):=p_{\sigma_1}(y)$ for all
$\sigma\in\Sigma^+$. Let $(e_1,...,e_n)^*$ denote a path. Then, by taking the natural version of $Z_n^{xy}$ and using $\log z\leq z-1$ for all $z>0$,
\begin{eqnarray*}
 \log Z_n^{xy}&=& \sum_{(e_1,...,e_n)^*}\log\frac{p_1^x...p_n^x}{p^y_1...p^y_n}1_{ _1[e_1,...,e_n]}\\
         &=& \sum_{(e_1,...,e_n)^*}\sum\limits_{i=1}^n\log\frac{p_i^x}{p_i^y}1_{ _1[e_1,...,e_n]}\\
         &\leq& \sum_{(e_1,...,e_n)^*}\sum\limits_{i=1}^n\frac{p_i^x-p_i^y}{p_i^y}1_{ _1[e_1,...,e_n]}\ \ \ \ P_x\mbox{-a.e.}.
\end{eqnarray*}

Now, observe that
\[\frac{p_i^x-p_i^y}{p_i^y}=\frac{p_i^x-p_i^y}{p_i^x}+\frac{(p_i^x-p_i^y)^2}{p_i^yp_i^x}.\]
Therefore,
\begin{equation}\label{JOE}
\log Z_n^{xy}\leq Y_n+X_n\ \ \ \ P_x\mbox{-a.e.}
\end{equation}
where
\[Y_n:=\sum\limits_{i=1}^n\sum_{(e_1,...,e_n)^*}\frac{p_i^x-p_i^y}{p_i^x}1_{ _1[e_1,...,e_n]}\]
and
\[X_n:=\sum\limits_{i=1}^n\sum_{(e_1,...,e_n)^*}\frac{(p_i^x-p_i^y)^2}{p_i^xp_i^y} 1_{ _1[e_1,...,e_n]}.\]
Furthermore, observe that
\[Y_{n+1}-Y_n=\sum_{(e_1,...,e_{n+1})^*}\frac{p_{n+1}^x-p_{n+1}^y}{p_{n+1}^x}1_{ _1[e_1,...,e_{n+1}]}\ \ \ \mbox{ for all }n\geq 1,\]
and, for every path $(e_1,...,e_{n})^*$,
\begin{eqnarray*}
&&\int\limits_{ _1[e_1,...,e_n]}(Y_{n+1}-Y_n)\ dP_x\\
&=&\sum\limits_{e_{n+1},i(e_{n+1})=t(e_{n})}\frac{p_{e_{n+1}}(w_{e_n}\circ...\circ
w_{e_1}x)-p_{e_{n+1}}(w_{e_n}\circ...\circ w_{e_1}y)}
{p_{e_{n+1}}(w_{e_n}\circ...\circ w_{e_1}x)}\\
&& \;\;\;\;\times p_{e_1}(x)...p_{e_{n}}(w_{e_{n-1}}\circ...\circ w_{e_1}x)p_{e_{n+1}}(w_{e_n}\circ...\circ w_{e_1}x)\\
&=&\sum\limits_{e_{n+1},i(e_{n+1})=t(e_{n})}(p_{e_{n+1}}(w_{e_n}\circ...\circ
w_{e_1}x)-p_{e_{n+1}}(w_{e_n}\circ...\circ w_{e_1}y))\\
&& \;\;\;\;\times p_{e_1}(x)...p_{e_{n}}(w_{e_{n-1}}\circ...\circ w_{e_1}x)\\
&=&0.
\end{eqnarray*}
Hence, $(Y_n,\mathcal{B}_n)_{n\in\mathbb{N}}$ is a $P_x$-martingale.
Therefore, $Y_n-Y_{n-1}$, $Y_{n-1}-Y_{n-2}$,..., $Y_2-Y_1$, $Y_1$
are orthogonal in $\mathcal{L}^2(P_x)$. By the Pythagoras equality,
this implies that
\begin{eqnarray*}
 \int Y^2_n dP_x&=&\int\left(\sum\limits_{i=2}^n(Y_i-Y_{i-1})+Y_1\right)^2 dP_x\\
  &=& \sum\limits_{i=2}^n\int(Y_i-Y_{i-1})^2\ dP_x+\int{Y_1}^2 dP_x\\
  &=&  \sum\limits_{i=1}^n\int\frac{(p_i^x-p^y_{i})^2}{(p_i^x)^2} dP_x.
\end{eqnarray*}

Therefore, 
\begin{eqnarray*}
 \int {Y_n}^2\ dP_x&\leq&\sum\limits_{i=1}^n\int\limits_{A_{xy}^{(\beta (i-1))}}\frac{1}{(p_i^x)^2}dP_x(\sigma)+\sum\limits_{i=1}^n\int\frac{ \phi^2_{\sigma_i}\left(a^{\beta (i-1)}d(x,y)\right) }{(p_i^x)^2}dP_x(\sigma)\\
&\leq&\sum\limits_{i=1}^n\sum\limits_{j\in N}\sum\limits_{\L_1[e_1,...,e_{i-1}]\subset A_{xy}^{(\beta (i-1))}, t(e_{i-1})=j}P_x\left(\L_1[e_1,...,e_{i-1}]\right)\ell_j\\
&&+\sum\limits_{i=1}^n\sum\limits_{j\in N}\phi^{(\beta (i-1))}_{xyj}{U^*}^{i-1}\delta_x(K_j).\\
\end{eqnarray*}
for all $n\in\mathbb{N}$. Hence, 
\[ \int {Y_n}^2\ dP_x\leq B_{\beta xy}\]
for all $n\in\mathbb{N}$.
The same upper bound holds for  $(\int {X_n}\ dP_x)_{n\in\N}$. Hence, by \eqref{JOE},
 \[H_n(P_x|P_y)=\int \log Z^{xy}_ndP_x\leq \left(\int{Y_n}^2dP_x\right)^\frac{1}{2}+\int {X_n}\ dP_x\leq \left(B_{\beta xy}\right)^\frac{1}{2}+B_{\beta xy}.\]
for all $n\in\N$.
\hfill$\Box$

\begin{Definition}
   Suppose $\M$ is positive and contractive with a contraction rate $0<a<1$.  Let $x,y\in K_i$ for some $i\in N$. For $0\leq \beta\leq 1$ and $k\geq 0$, set  
\[\hat\ell_{xk}:=\sup\{\ell_j|\ j\mbox{ is accessible from $i$ by a path of length }k\},\]
\[\hat\phi^{(\beta k)}_{xy}:=\sup\left\{\left.\phi^{(\beta k)}_{xyj}\right|\ j\mbox{ is accessible from $i$ by a path of length }k\right\}\]
and, for $t\geq 0$,
\[\phi_\M(t):=\sup\limits_{e\in E}\phi_e(t).\]
\end{Definition}

\begin{prop}\label{reub}
     Suppose $\M$ is positive and contractive with a contraction rate $0<a<1$. Let $x,y\in K_i$ for some $i\in N$ and $0\leq \beta\leq 1$. Then the following holds true.\\
(i)  \[B_{\beta xy}\leq\sum\limits_{k=0}^\infty\left(\hat\ell_{xk}a^{(1-\beta)k}+\hat\phi^{(\beta k)}_{xy}\right).\]
(ii) \[B_{\beta xy}\leq\sum\limits_{k=0}^\infty\hat\ell_{xk}\left(a^{(1-\beta)k}+\phi^2_\M \left(a^{\beta k}d(x,y)\right)\right).\]
(iii) If each $w_e|_{K_{i(e)}}$ is  contractive with a contraction rate $0<a<1$, then 
\[B_{\beta xy}\leq\sum\limits_{k=0}^\infty\hat\ell_{xk}\phi^2_\M \left(a^{\beta k}d(x,y)\right).\]
\end{prop}
{\it Proof.} 
By the contraction on average condition,
\[\int d(w_{\sigma_{k}}\circ...\circ w_{\sigma_1}(x),w_{\sigma_{k}}\circ...\circ w_{\sigma_1}(y))\ dP_x\leq a^kd(y,x)\mbox{ for all }
k\geq 0.\] 
Hence,
\[P_x\left(A_{xy}^{(\beta k)}\right)\leq a^{(1-\beta)k}\mbox{ for all }k\geq 0.\]
If all $w_e|_{K_{i(e)}}$ are contractive with a contraction rate $0<a<1$, then $A_{xy}^{(\beta k)}$ is empty for all $k\geq 0$. Thus, the assertions follow immediately from the definition of $B_{\beta xy}$.
  \hfill$\Box$

\begin{Example}\label{se}
   Let $\mathcal{D}_R:=([0,1], w'_e,p'_e)_{e\in \{0,1\}}$ where
\begin{eqnarray*}
    w'_0(z):=\frac{1}{2}z,&\ & w'_1(z):=\frac{1}{2}+\frac{1}{2}z,\\
    p'_0(z):=z,&\ &p'_1(z):=1-z
\end{eqnarray*}
for all $z\in [0,1]$. Let $\M:=(K_{i(e)},w_e, p_e)_{e\in E}$ be the Markov system resulting from $\mathcal{D}_R$ through the restriction of the maps and the probability functions on the atoms of the following Markov partition. Set $K_0:=\{0\}$,  $K_\infty:=\{1\}$, $K_j:=(1-1/2^j,1-1/2^{j+1}]$ for all $j\in\N$ and $K_j:=(1/2^{|j|+1},1/2^{|j|}]$ for all $j\in\Z\setminus(\N\cup\{0\})$. Note that, since $1/p'_0(z)+1/(1-p'_1(z))=1/(z(1-z))$ for all $z\in]0,1[$, $\ell_j\leq 2^{|j|+2}$ for all $j\in\Z\setminus\{0\}$. 

Now, let $i_0\in\Z\setminus\{0\}$ and $x,y\in K_{i_0}$. From the directed graph associated with $\M$, one sees that $\hat\ell_{xk}\leq2^{|i_0|+k+2}$ for all $k\geq 0$. Hence, by Proposition \ref{reub} (iii), for $\beta=1$,
\[B_{1 xy}\leq2^{|i_0|+2}|x-y|^2\sum\limits_{k=0}^\infty 2^k \left(\frac{1}{2}\right)^{2 k}=2^{|i_0|+3}|x-y|^2,\]
and therefore, by Lemma \ref{acl},
\[H(P_x|P_y)\leq2^{\frac{|i_0|+3}{2}}|x-y|+2^{|i_0|+3}|x-y|^2.\]
Thus, by Theorem \ref{klfc}, $\M$ satisfies Condition \ref{acC}. Furthermore, one sees the same way as in Example 3 in \cite{Wer11} that $\M$ is non-degenerate, as
 \[R1=...+1_{\left\{\frac{1}{8}\right\}}+1_{\left\{\frac{1}{4}\right\}}+1_{\left\{\frac{1}{2}\right\}}+1_{\left\{1-\frac{1}{4}\right\}}+1_{\left\{1-\frac{1}{8}\right\}}+...\] 
and 
 \[R^21=...+\frac{1}{4}1_{\left\{\frac{1}{8}\right\}}+\frac{1}{2}1_{\left\{\frac{1}{4}\right\}}+\frac{1}{2}1_{\left\{1-\frac{1}{4}\right\}}+\frac{1}{4}1_{\left\{1-\frac{1}{8}\right\}}+....\] 
Also, obviously, $\M$  satisfies  \eqref{cgc} for any choice of $x_j\in K_j$ for all $j\in N$.
Hence, $\M$ satisfies conditions of Theorem \ref{esgc}. Thus, in particular, it has a unique invariant Borel probability measure, and therefore, the same is true for $\mathcal{D}_R$.

Clearly,  $\M$ also satisfies the conditions of Corollary \ref{Ec}. We show that, in this case, the functions $f_e(x):=-\log p_e(x)$ if $x\in K_{i(e)}$ and $f(x):=0$ otherwise also satisfy the conditions of Theorem \ref{esgc} (ii) (a). By Remark \ref{ekcr}, $\sum_{e\in\{0,1\}}p_e(\log p_e)^2\leq 2(\log 2)^2$. Define the function  $h:=-\sum_{e\in\{0,1\}}p_e\log p_e$. Then 
\[\frac{dh}{dx}=-\log\frac{x}{1-x}\ \ \ \mbox{ for all }x\in(0,1).\]
A simple computation shows that
\[\left|\log\frac{x}{1-x}\right|\leq (|j|+1)\log2\ \ \ \mbox{ for all }x\in K_j\mbox{ and }j\in\Z\setminus\{0\}.\]
Hence
\[\Delta_jh(t)\leq t(|j|+1)\log 2\ \ \ \mbox{ for all }t>0,\]
and therefore,
\[\Delta^{(k)}_{i_0}h\left(2^{-\alpha k}d(x,y)\right)\leq 2^{-\alpha k}d(x,y)(|i_0|+k+3)\log 2\ \ \ \mbox{ for all }k\in\N\mbox{ and }0<\alpha<1.\]
Hence, $h$ is $(\M,\alpha)$-uniformly continuous for all $0<\alpha<1$. 

Now, consider the Markov system $\M_2$ associated with $D_R$ which results from the partition $(K^2_j)_{j\in\Z\cup\{\infty\}}$ where $K^2_0:=\{0\}$,  $K^2_\infty:=\{1\}$, $K^2_j:=[1-1/2^j,1-1/2^{j+1})$ for all $j\in\N$ and $K^2_j:=[1/2^{|j|+1},1/2^{|j|})$ for all $j\in\Z\setminus(\N\cup\{0\})$. Obviously, the same way as above, $\M_2$ also satisfies Condition \ref{acC}. Now, note that $\M$ and $\M_2$ are connecting refinements of $\M_0:=(K^0_{i(e)}, p^0_e,w^0_e)_{e\in\{a,b,c,d\}}$ given by the restrictions of the maps and the probability functions of $D_R$ on $K^0_0:=\{0\}$, $K^0_{1}:=(0,1)$ and $K^0_2:=\{1\}$.
Thus, by Lemma \ref{crl}, $\M_0$  satisfies the conditions of Theorem \ref{esgc} and Corollary \ref{Ec} also.
\end{Example}

\subsection*{Acknowledgements} 
The author would like to thank an anonymous reviewer for the Annales de l'Institut Henri Poincar\'{e} for suggestions on improvements to the text of this article, an anonymous reviewer for the Journal of Modern Dynamics for numerous corrections of misprints and grammar and suggestions on further improvements to the text of the article,  Boris M. Gurevich for the invitations to give several talks on the subject at the Ergodic Theory and Statistical Mechanics Seminar at the Lomonosov Moscow State University and also other participants of the seminar for valuable comments and questions which helped to improve the article.

\end{document}